\documentclass[fleqn]{mat01}
\usepackage{times,mathtimy,amssymb,latexsym}

\newcommand{\ord}{a}
\newcommand{\ga}{\gamma}
\newcommand{\eps}{\epsilon}
\newcommand{\al}{\alpha}

\newcommand{\la}{\lambda}
\newcommand{\D}{\mathfrak{F}}

\newcommand{\RR}{\mathbb{R}}
\newcommand{\lap}{\Delta}
\newcommand{\LL}{\mathfrak{L}}
\newcommand{\hf}{h}
\newcommand{\NN}{\mathbb{N}_0}

\newcommand{\CC}{C_c^\infty}

\begin{document}

\newtheorem{theor}{\bf Theorem}
\newtheorem{lem}{\it Lemma}
\newtheorem{propo}{\rm PROPOSITION}
\newtheorem{coro}{\rm COROLLARY}
\newtheorem{defini}{\rm DEFINITION}
\newtheorem{rem}{\it Remark}
\newtheorem{poft}{\it Proof of Theorem}

\setcounter{page}{337} \firstpage{337}

\title{Sobolev spaces associated to the harmonic oscillator}

\markboth{B~Bongioanni and J~L~Torrea}{Sobolev spaces associated
to the harmonic oscillator}

\author{B~BONGIOANNI$^{1}$ and J~L~TORREA$^{2}$}

\address{$^{1}$Departamento de Matem\'atica, Facultad de Ingenier\'{\i}a
Q\'{\i}mica, Universidad Nacional del Litoral, and Instituto de
Matem\'atica Aplicada del Litoral, Santa Fe (3000), Argentina\\
\noindent $^{2}$Departamento de Matem\'atica, Facultad de
Ciencias, Universidad Aut\'onoma de Madrid, Spain\\
\noindent E-mail: bbongio@math.unl.edu.ar; joseluis.torrea@uam.es}

\volume{116}

\mon{August}

\parts{3}

\pubyear{2006}

\Date{MS received 27 September 2005}

\begin{abstract}
We define the Hermite--Sobolev spaces naturally associated to the
harmonic oscillator $H= -\Delta +|x|^2.$ Structural properties,
relations with the classical Sobolev spaces, boundedness of
operators and almost everywhere convergence of solutions of the
Schr\"odinger equation are also considered.
\end{abstract}


\maketitle

\section{Introduction}

We consider the second-order differential operator
\begin{equation}\label{Hop}
H=-\Delta +|x|^2,\quad x \in \RR^d.
\end{equation}
This operator is self-adjoint on the set of infinitely
differentiable functions with compact support $\CC$, and it can be
factorized as
\begin{equation}\label{Hfact}
H= \frac{1}{2} \sum_{j=1}^d A_jA_{-j}+ A_{-j}A_{j},
\end{equation}
where
\begin{equation*}
A_j = \frac{\partial}{\partial x_j} + x_j\quad \text{and}\quad
A_{-j} = A_j^* = -\frac{\partial}{\partial x_j} + x_j,\quad
\text{$1\leq j\leq d$.}
\end{equation*}

In the last few years several authors have been concerned with the
harmonic analysis associated to the operator $H$ (see for instance
\cite{5,10,12}). In this analysis the operators $A_j$ play the
role of the partial derivative operators ${\partial}/{\partial
x_j}$ in the classical Euclidean case. Hence it seems natural to
study the spaces of functions in $L^p(\RR^d)$ whose {\em
derivatives} also belong to $L^p(\RR^d).$ Following this idea, we
introduce the Hermite--Sobolev spaces $W^{k,p}$
(Definition~\ref{defi:sob}). These spaces are Banach spaces and
the set of  linear combinations of Hermite functions is dense in
any of them (Proposition~\ref{propo:densidad}). The spaces
$W^{k,p}$ were previously studied in \cite{14} for $p=2$ and in
\cite{7} for $p\not=2$.

Once we have the {\em Laplacian} $H$, it is also natural to
consider the {\em potential spaces} ${\LL}_a^p =
H^{-a/2}(L^p(\RR^d))$ (Definition~\ref{defi:esppot}). In other
words, the range of the {\em Hermite fractional integral} operator
$H^{-a/2}$  in $L^p(\RR^d).$ In order to have a satisfactory
description of these potential spaces we need a sharp analysis of
the operator $H^{-a},$ for $a>0.$ Such analysis is contained in
Proposition~\ref{propo:KerHal} and Lemma~\ref{lem:xaLasobre2}. It
turns out that, for $k \in \mathbb{N}$,  the spaces ${\LL}_k^p$
and $W^{k,p}$ coincide (see Theorem~\ref{teo:Wkp}). The proof of
this theorem uses the boundedness in $L^p(\RR^d)$ of the Riesz
transforms naturally associated to $H.$ These Riesz transforms
were introduced by Thangavelu in \cite{12}, and some of their
boundedness properties can be found in \cite{5} and \cite{10}.

Observe that, in some sense, Theorem~\ref{teo:Wkp} allows us to
say that the spaces ${\LL}_\ord^p$ are the spaces of functions in
$L^p(\RR^d)$ for which their {\em derivatives of order $\ord$}
also belong to $L^p(\RR^d).$

Once we have a satisfactory definition of Hermite--Sobolev (or
Hermite {\it potential}) spaces and hence of {\it fractional}
derivatives, we study their relationship with the corresponding
classical Euclidean spaces. We show in Theorem~\ref{teo:HSvsS}
that although the Hermite--Sobolev spaces coincide locally with
the Euclidean Sobolev spaces, they are in fact different. In
\S\ref{sec:boud_op_LL}, we show that the Hermite--Riesz transforms
are bounded on the Hermite--Sobolev spaces while the classical
Hilbert transform is not bounded on these spaces. From the careful
analysis of the kernel of $H^{-\ord},$ we also obtain certain
inequalities of {\it Poincar\'e type} for the {\it derivatives}
$A_j$ in Theorem~\ref{teo:poinca}. Finally, in
\S\ref{sec:schrodinger} we give an application to the almost
everywhere convergence of  the solution of the Schr\"odinger
equation~\eqref{scho_eq} to the initial data.

Our work was heavily inspired in the paper by Thangavelu
\cite{14}, where the spaces ${\LL}_\ord^2$ were defined. These
spaces were also considered in \cite{6}.

As we said above, in order to develop this work, some nontrivial
estimates of the Hermite fractional integral operator
(Definition~\ref{defi:Ha}) were needed. However, it is not the aim
of this paper to make an exhaustive study of this operator. Hence,
other natural questions like weak and strong boundedness in the
extreme points or $BMO$-type boundedness of the operator $H^{-a}$
are left aside and they will be the motivation of a forthcoming
paper.

\section{Hermite--Sobolev spaces}

Let $n\in \NN=\mathbb{N}\cup \{0\}$ and consider the Hermite
function of order $n$,
\begin{equation*}
\hf_n(t)= \frac{(-1)^n}{(2^n n! \pi^{1/2})^{1/2}}\, H_n(t)\,
\hbox{e}^{-t^2/2},\quad t\in\RR,
\end{equation*}
where $H_n$ denotes the Hermite polynomial of degree $n$ (see
\cite{12}). Given a multi-index $\al=(\al_j)_{j=1}^{d}\in \NN^d$,
we consider the Hermite function, $\hf_\al$, as
\begin{equation*}
\hf_\al (x)=\prod_{j=1}^d \hf_{\al_j}(x_j),\quad x=(x_1,\dots,x_d)
\in\RR^d.
\end{equation*}
These functions are eigenvectors of the Hermite operator $H$
defined in \eqref{Hop}. In fact
\begin{equation*}
H \hf_\al = (2|\al|+d)\, \hf_\al,
\end{equation*}
where $|\al|=\sum_{j=1}^d \al_j$. Moreover, for $1\leq j\leq d$,
\begin{equation*}
A_j\hf_\al = \sqrt{2\al_j}\, \hf_{\al-e_j},\quad A_{-j}\hf_\al =
\sqrt{2(\al_j+1)}\, \hf_{\al+e_j},
\end{equation*}
where $e_j$ is the $j$th coordinate vector in $\NN^d$. The
operators $A_j$ and $A_{-j}$ are called \emph{annihilation} and
\emph{creation} operators respectively.

\begin{defini}\label{defi:sob}$\left.\right.$\vspace{.5pc}

\noindent {\rm Given $p\in (1,\infty)$ and $k\in \mathbb{N}$, we
define the Hermite--Sobolev space of order $k$, denoted by
$W^{k,p}$, as the set of functions $f\in L^p(\RR^d)$ such that
\begin{equation*}
A_{j_1}\cdots A_{j_m}f \in L^p(\RR^d),\quad
1\leq|j_1|,\dots,|j_m|\leq d, \quad 1\leq m\leq k,
\end{equation*}
with the norm
\begin{equation*}
\|f\|_{W^{k,p}} = \sum_{1\leq |j_1|,\dots,|j_m| \leq d ,\ 1\leq
m\leq k} \|A_{j_1}\cdots A_{j_m}f\|_p + \|f\|_p\text{.}
\end{equation*}
}
\end{defini}

We will show that the set of finite linear combinations of Hermite
functions, denoted by $\D$, is dense in the Hermite--Sobolev
spaces. We shall need the following lemmas. Their proofs may be
found in \cite{10} and \cite{12}, respectively.

\begin{lem}\label{lem:StemTorr}
Let $m\in\NN$ and $f\in\CC$. There exists a constant $C_{m,f}>0$
such that
\begin{equation*}
|\langle f,\hf_\al \rangle| \leq C_{m,f}\,(|\al|+1)^m\text{,}\ \ \
\al\in\NN^d\text{.}
\end{equation*}
\end{lem}

\begin{lem}\label{lem:thangavelu}
As $n\rightarrow\infty$ the Hermite functions satisfy the
estimates
\begin{enumerate}
\renewcommand\labelenumi{\rm (\roman{enumi})}
\leftskip .4pc
\item $\|\hf_n\|_p \sim n^{\frac{1}{2p}-\frac{1}{4}},\quad
1\leq p < 4${\rm ,}

\item $\|\hf_n\|_p \sim n^{-\frac{1}{8}}\log(n),\quad p =
4${\rm ,}

\item $\|\hf_n\|_p \sim n^{\frac{1}{6p}-\frac{1}{12}},\quad
4< p \leq \infty$.
\end{enumerate}
\end{lem}

\begin{propo}\label{propo:densidad}$\left.\right.$\vspace{.5pc}

\noindent Let $p$ be in the range $1<p<\infty$ and $k\in
\mathbb{N}$. The set $W^k_p$ is a Banach space. Moreover{\rm ,}
the sets $\D$ and $\CC$ are dense in $W^{k,p}$.
\end{propo}

\begin{proof}
To see that $W^{1,p}$ is complete, observe that if $\{f_n\}_{n\ge
1}$ is a Cauchy sequence in $W^{1,p}$, then
\begin{equation}\label{suces}
\left\{ \frac{\partial f_n}{\partial x_j}\right\}_{n\ge 1}\quad
\hbox{and}\quad \{x_jf_n\}_{n\ge 1},\quad 1\leq j\leq d\text{,}
\end{equation}
are Cauchy sequences in $L^p(\RR^d)$. If we call $f$ the limit in
$L^p(\RR^d)$ of $\{f_n\}_{n\ge 1}$, it is easy to see that
${\partial f}/{\partial x_j}$ and $x_jf$ are respectively the
limits in $L^p(\RR^d)$ of \eqref{suces} (see \cite{8}, p.~122).

Now we shall see that $\CC$ is a dense set in $W^{1,p}$. Let
$\psi$ be a function in $\CC$ such that $\int_{\RR^d} \psi = 1$.
For every $\eps>0$, consider
\begin{equation*}
\psi_\eps(x) = \frac{1}{\eps^d}\,
\psi\!\left(\frac{x}{\eps}\right).
\end{equation*}
Given $f$ in $W^{1,p}$, define $f_\eps = f * \psi_\eps$. Following
the ideas in p.~123 of \cite{8}, we have
\begin{equation*}
\|f-f_\eps\|_{p} \rightarrow 0\quad \hbox{and}\quad
\left\|\frac{\partial}{\partial x_j}f_\eps -
\frac{\partial}{\partial x_j}f \right\|_p \rightarrow 0, \quad
1\leq j \leq d.
\end{equation*}
On the other hand, for $1\leq j\leq d$, we call
$\phi^j(x)=x_j\,\psi(x)$ and  $\phi^j_\eps(x) = \frac{1}{\eps^d}\,
\phi^j\!\left(\frac{x}{\eps}\right)$, then the function $\phi^j\in
\CC$, and for $\eps>0$,
\begin{align}\label{xfpsi_xfe}
\|(x_jf)*\psi_\eps - x_jf_\eps\|_p^p &= \int_{\RR^d}
\left|\int_{\RR^d} f(y)(y_j-x_j)
\frac{1}{\eps^d}\psi\left(\frac{x-y}{\eps}\right)\,\hbox{d}y\right|^p
\hbox{d}x\nonumber\\[.5pc]
&= \eps^p \int_{\RR^d} \left|\int_{\RR^d} f(y) \frac{1}{\eps^d}
\phi^j\left(\frac{x-y}{\eps}\right)\,\hbox{d}y\right|^p \hbox{d}x\nonumber\\[.5pc]
&= \eps^p\ \|f * \phi^j_\eps\|_p^p.
\end{align}
Moreover, since $x_jf$ belongs to $L^p(\RR^d)$,
\begin{equation}\label{xfpsi_xf}
\|(x_jf)*\psi_\eps - x_jf\|_p \rightarrow 0.
\end{equation}
Equations~\eqref{xfpsi_xfe} and \eqref{xfpsi_xf} imply $x_jf_\eps
\ \rightarrow\ x_jf$ in $L^p(\RR^d)$. Therefore, we conclude that
$f_\eps\in W^{1,p}$ and $f_\eps\rightarrow f$ in the
$W^{1,p}$-norm.

The functions $f_\eps$ do not necessarily have compact support,
but they can be modified as in the classical case (see \cite{8},
p.~123).

It remains to prove that any function in $\CC$ can be approximated
in the $W^{1,p}$-norm by a function in $\D$. In fact, we will show
that any $f\in\CC$ is the limit, in the $W^{p,1}$-norm, of a
subsequence of the partial sums
\begin{equation*}
S_N f = \sum_{|\al|\leq N} \langle f,\hf_{\al}\rangle\,
\hf_\al\text{.}
\end{equation*}
In \cite{10}, it is proved that there exists a subsequence of the
previous sequence converging to $f$ in the $L^p$-norm. Hence, it
is enough to show that there exists a subsequence of
$\{A_j(S_N(f))\}_{N\ge 1}=\{S_N(A_jf)\}_{N\ge 1}$ converging to
$A_jf$ in the $L^p$-norm, for every $j$ with $1\leq |j|\leq d$.

Let us fix $j$ to be in $1\leq j\leq d$ (the case $-d\leq j\leq
-1$ is similar). The sequence $\{S_N(A_jf)\}_{N\ge 1}$ converges
to $A_j f$ in the $L^2$-norm. Hence we can take a subsequence
$\{S_{N_k}(A_jf)\}_{k\ge 1}$ converging to $A_j f$ almost
everywhere. Since
\begin{align*}
S_N(A_jf) & = \sum_{|\al|\leq N} \langle A_j f,\hf_\al \rangle
\hf_\al = \sum_{|\al|\leq N} \langle f, A_{-j} \hf_\al\rangle
\hf_\al\\[.5pc]
& = \sum_{|\al|\leq N} \sqrt{2(\al_{-j}+1)} \langle
f,\hf_{\al+e_j}\rangle \hf_\al,
\end{align*}
by Lemma~\ref{lem:StemTorr} and Hölder's inequality,
\begin{align*}
|S_N(A_jf)|^p &\leq C \left(\sum_{|\al|\leq N}
\sqrt{2(\al_{-j}+1)} (|\al|+2)^{-M}|\hf_\al|\right)^p\\[.5pc]
&\leq C \left(\sum_{|\al|\leq N} (|\al|+1)^{-M+1}
|\hf_\al|\right)^p
\end{align*}
\begin{align*}
&\leq C \left(\sum_{|\al|\leq N}
(|\al|+1)^{-M+1}\right)^{p/p'}\!\!\sum_{|\al|\leq N}
(|\al|+1)^{-M+1} |\hf_\al|^p\\[.5pc]
&\leq C\, \sum_{\al} (|\al|+1)^{-M+1} |\hf_\al|^p.
\end{align*}
From Lemma~\ref{lem:thangavelu}, it is easy to see that the
function $\sum_{\al} (|\al|+1)^{-M+1} |\hf_\al|^p$ belongs to
$L^1(\RR^d,\hbox{d}x)$, and the dominated convergence theorem
gives that $\{S_{N_k}(A_jf)\}_{k\ge 1}$ converges to $A_j f$ in
the $L^p$-norm. Now we can repeat the argument for every $j$,
taking a subsequence of the previous subsequence in each step.

For $k>1$ we leave the details to the reader.\hfill $\Box$
\end{proof}

\section{Fractional integral}

With the ideas in \cite{8} and \cite{10}, we introduce the
following operator.

\begin{defini}\label{defi:Ha}$\left.\right.$\vspace{.5pc}

\noindent {\rm  Given  $\ord> 0$, we define for $f\in\D$, the
operator
\begin{equation}\label{Haf}
H^{-\ord}f(x) = \frac{1}{\Gamma(\ord)} \int_0^\infty
\hbox{e}^{-tH}f(x)\,t^{\ord}\,\frac{\hbox{d}t}{t},\quad x\in\RR^d,
\end{equation}
where $\{\hbox{e}^{-tH}\}_{t\ge 0}$ is the heat semi-group
associated to $H$.}
\end{defini}

\begin{rem}{\rm
If $\ord>0$ and $\al\in \NN^d$, by using the $\Gamma$ function and
the fact
\begin{equation*}
\hbox{e}^{-tH} h_\al = \hbox{e}^{-t(2|\alpha|+d)}h_\al,
\end{equation*}
we have
\begin{equation*}
H^{-\ord}\hf_\al(x) = \frac{1}{\Gamma(\ord)} \int_0^\infty
\hbox{e}^{-tH}\hf_\al(x)\,t^{\ord}\,\frac{\hbox{d}t}{t}=
(2|\al|+d)^{-\ord}\hf_\al(x),\quad x\in\RR^d.
\end{equation*}}
\end{rem}

\begin{propo}\label{propo:KerHal}$\left.\right.$\vspace{.5pc}

\noindent The operator $H^{-\ord}$ has integral representation
\begin{equation}\label{Hal}
H^{-\ord}f(x) =  \int_{\RR^d}\! K_\ord(x,y)\,f(y)\,\hbox{\rm
d}y,\quad x\in\RR^d,
\end{equation}
for all $f\in \D$. Moreover{\rm ,} there exist $\Phi_\ord$ in
$L^1(\RR^d)$ and a constant $C$ such that
\begin{equation}\label{acotPhi}
K_\ord(x,y) \leq C\,\Phi_\ord(x-y),\quad \hbox{for all} \
x,y\in\RR^d.
\end{equation}
\end{propo}

\begin{proof}
If $f\in\D$, then for $x\in\RR^d$,
\begin{align*}
H^{-\ord}f(x) &= \frac{1}{\Gamma(\ord)} \int_0^\infty
\hbox{e}^{-tH}f(x)\,t^{\ord}\,\frac{\hbox{d}t}{t}\\[.5pc]
&= \frac{1}{\Gamma(\ord)} \int_0^\infty \int_{\RR^d}
G_t(x,y)\,f(y)\,\hbox{d}y\,t^{\ord}\,\frac{\hbox{d}t}{t},
\end{align*}
where
\begin{equation*}
G_t(x,y) = (2\pi \sinh 2t)^{-d/2} \hbox{e}^{-\frac{1}{2} |x-y|^2
\coth 2t - x\cdot y\tanh t},
\end{equation*}
for $x,y\in\RR^d$, (see \cite{10}). Therefore, if we show that for
some constant $C$,
\begin{equation*}
\frac{1}{\Gamma(\ord)}\int_0^\infty G_t(x,y)\,
t^\ord\,\frac{\hbox{d}t}{t}\leq C\,\Phi_\ord(x-y),\quad \hbox{for
all} \ x,y\in\RR^d,
\end{equation*}
where $\Phi_\ord\in L^1(\RR^d)$, then by Fubini's theorem,
\begin{equation*}
H^{-\ord}f(x) =  \int_{\RR^d}\! K_\ord(x,y)\,f(y)\,\hbox{d}y,\quad
x\in\RR^d,
\end{equation*}
with
\begin{equation}\label{Kker}
K_\ord(x,y) = \frac{1}{\Gamma(\ord)}\int_{0}^\infty
G_t(x,y)\,t^{\ord}\frac{\hbox{d}t}{t}.
\end{equation}
We perform the change of variables
\begin{equation*}
t=\frac{1}{2} \log\left(\frac{1+s}{1-s}\right),
\end{equation*}
then
\begin{equation}\label{Kkercambio}
K_\ord(x,y)= \frac{1}{\Gamma(\ord)(4\pi)^{d/2}2^{\ord-1}} \int_0^1
\zeta_\ord (s)\, \hbox{e}^{-\frac{1}{4}
\left(s|x+y|^2+\frac{1}{s}|x-y|^2\right)}\, \frac{\hbox{d}s}{s},
\end{equation}
where
\begin{equation*}
\zeta_\ord (s)=\left(\frac{1-s^2}{s} \right)^{\frac{d}{2}-1}\!\!
\log\left(\frac{1+s}{1-s}\right)^{\ord-1}.
\end{equation*}

We split $K_\ord$ as $K_\ord = K_{\ord,0} + K_{\ord,1}$, where
\begin{equation}\label{K0}
K_{\ord,0}(x,y) = \frac{1}{\Gamma(\ord)(4\pi)^{d/2}
2^{\ord-1}}\int_0^{1/2} \zeta_\ord (s)\,
\hbox{e}^{-\frac{1}{4}\left(s|x+y|^2+\frac{1}{s}|x-y|^2\right)}\,
\frac{\hbox{d}s}{s}.
\end{equation}
Since the integral
\begin{equation*}
\int_{1/2}^1 \zeta_\ord (s)\, \frac{\hbox{d}s}{s}\ < \infty,
\end{equation*}
we have
\begin{align}\label{K1}
K_{\ord,1}(x,y) &=
\frac{1}{\Gamma(\ord)(4\pi)^{d/2}2^{\ord-1}}\int_{1/2}^1
\zeta_\ord (s)\, \hbox{e}^{-\frac{1}{4} \left(s|x+y|^2+
\frac{1}{s}|x-y|^2\right)}\, \frac{\hbox{d}s}{s}\nonumber\\[.4pc]
&\leq C\,\hbox{e}^{-\frac{1}{4}|x-y|^2}\,
\hbox{e}^{-\frac{1}{8}|x+y|^2}.
\end{align}

It is easy to see that there exists a constant $C_1$ which depends
on $d$ and $\ord$ such that
\begin{equation}\label{Z_0y1med}
\frac{s^{-\frac{d}{2}+\ord}}{C_1} \leq \zeta_\ord (s) \leq C_1
s^{-\frac{d}{2}+\ord}\quad \hbox{for} \ 0<s<1/2.
\end{equation}
Therefore
\begin{equation*}
K_{\ord,0}(x,y) \leq  C_1 \int_0^{1/2} s^{-\frac{d}{2}+\ord}\,
\hbox{e}^{-\frac{1}{4s}|x-y|^2}\, \frac{\hbox{d}s}{s}.
\end{equation*}

If $|x-y|\ge 1$, the last expression is bounded up to a constant
by
\begin{equation*}
\hbox{e}^{-\frac{1}{4}|x-y|^2} \int_0^{1/2}
s^{-\frac{d}{2}+\ord}\, \hbox{e}^{-\frac{1}{8s}}\,
\frac{\hbox{d}s}{s},
\end{equation*}
hence,
\begin{equation}\label{Kacot-glob}
K_{\ord,0}(x,y) \leq  C\,\hbox{e}^{-\frac{1}{4}|x-y|^2}.
\end{equation}

Let us study now the region  $|x-y|<1.$ By a change of variables,
\begin{equation*}
\int_0^{1/2} s^{-\frac{d}{2}+\ord}\,
\hbox{e}^{-\frac{1}{4s}|x-y|^2}\, \frac{\hbox{d}s}{s} =
\frac{1}{|x-y|^{d-2\ord}} \int_{2|x-y|^2}^\infty
s^{\frac{d}{2}-\ord}\, \hbox{e}^{-\frac{s}{4}}\,
\frac{\hbox{d}s}{s}.
\end{equation*}
In the case $\ord< (d/2)$, we get
\begin{equation}\label{Kacot-loc-alp}
K_{\ord,0}(x,y) \leq  \frac{C}{|x-y|^{d-2\ord}} \text{.}
\end{equation}
For the case $\ord = (d/2)$,
\begin{align*}
\int_{2|x-y|^2}^\infty \hbox{e}^{-\frac{s}{4}}\,
\frac{\hbox{d}s}{s} &= \int_{2|x-y|^2}^2 \frac{\hbox{d}s}{s} +
\int_2^\infty \hbox{e}^{-\frac{s}{4}}\,
\frac{\hbox{d}s}{s}\\[.5pc]
&= \log\left(\frac{1}{|x-y|^2}\right) +
\int_2^\infty \hbox{e}^{-\frac{s}{4}}\, \frac{\hbox{d}s}{s}\\[.5pc]
&\leq  2 \log\left(\frac{e}{|x-y|}\right) \left( 1 + \int_2^\infty
\hbox{e}^{-\frac{s}{4}}\, \frac{\hbox{d}s}{s}\right).
\end{align*}
Then,
\begin{equation}\label{Kacot-loc-allog}
K_{\ord,0}(x,y) \leq  C\,\log\left(\frac{e}{|x-y|}\right)\text{.}
\end{equation}
Finally, when $\ord> (d/2)$,
\begin{align*}
\int_{2|x-y|^2}^\infty s^{\frac{d}{2}-\ord}\,
\hbox{e}^{-\frac{s}{4}}\, \frac{\hbox{d}s}{s} & \leq
\int_{2|x-y|^2}^2 s^{\frac{d}{2}-\ord}\, \frac{\hbox{d}s}{s} +
\int_2^\infty s^{\frac{d}{2}-\ord}\,
\hbox{e}^{-\frac{s}{4}}\, \frac{\hbox{d}s}{s}\\[.5pc]
&\leq\frac{1}{2^{2\ord-d}\left(\ord-\frac{d}{2}\right)|x-y|^{2\ord-d}}
+ \int_2^\infty s^{\frac{d}{2}-\ord}\, \hbox{e}^{-\frac{s}{4}}\,
\frac{\hbox{d}s}{s}.
\end{align*}
Thus,
\begin{equation}\label{Kacot-loc-alg}
K_{\ord,0}(x,y) \leq C.
\end{equation}

Therefore, from \eqref{K1},
\eqref{Kacot-glob}--\eqref{Kacot-loc-alg}, if we define for
$x\in\RR^d$,
\begin{equation}\label{laPhia}
\Phi_\ord(x)= \begin{cases}
\frac{\chi_{\{|x|<1\}}(x)}{|x|^{d-2\ord}} +
\hbox{e}^{\frac{-|x|^2}{4}}\,
\chi_{\{|x|\ge1\}}(x), &\hbox{if} \ \ord<\frac{d}{2},\\[.5pc]
\log\left(\frac{e}{|x|}\right)\,\chi_{\{|x|<1\}}(x) +
\hbox{e}^{\frac{-|x|^2}{4}}\, \chi_{\{|x|\ge1\}}(x), &\hbox{if} \
\ord=\frac{d}{2},\\[.5pc]
\chi_{\{|x|<1\}}(x) + \hbox{e}^{\frac{-|x|^2}{4}}\,
\chi_{\{|x|\ge1\}}(x), &\hbox{if} \ \ord>\frac{d}{2},
\end{cases}
\end{equation}
we have proved \eqref{acotPhi}.\hfill $\Box$
\end{proof}

\begin{theor}[\!]\label{teo:fcoef}
The operator $H^{-\ord}$ defined by {\rm \eqref{Haf},} is
well-defined and bounded on $L^p(\RR^d)${\rm ,} for all
$p\in[1,\infty]$. Moreover{\rm ,} for all $f$ in $L^p(\RR^d)$ and
$\al\in\NN^d${\rm ,} we have
\begin{equation}\label{fcoef}
\int_{\RR^d} H^{-\ord}f\, \hf_\al = (2|\al|+d)^{-\ord}
\int_{\RR^d} f\,\hf_\al.
\end{equation}
\end{theor}

\begin{proof}
The boundedness of the operator $H^{-\ord}$ on $L^p(\RR^d)$ is due
to the fact that the kernel $K_\ord$ is bounded by an integrable
function. To see \eqref{fcoef}, let $\al\in\NN^d$. By
Proposition~\ref{propo:KerHal} and H\"older's inequality,
\begin{align*}
&\int_{\RR^d}\! \int_{\RR^d}\!
|K_\ord(x,y)|\,|f(y)|\,|\hf_\al(x)|\,\hbox{d}y\,\hbox{d}x\\[.4pc]
&\quad\, \leq C\,\int_{\RR^d}\! \int_{\RR^d}\!
|\Phi_\ord(x-y)|\,|f(y)|\,|\hf_\al(x)|\,\hbox{d}y\,\hbox{d}x\\[.4pc]
&\quad\, \leq C\, \|f\|_p \|\hf_\al\|_{p'},
\end{align*}
where the constant $C$ depends on $\al$. Therefore, by Fubini's
theorem,
\begin{align*}
\int_{\RR^d} H^{-\ord}f\, \hf_\al = \int_{\RR^d}\!
f\,H^{-\ord}\hf_\al = (2|\al|+d)^{-\ord} \int_{\RR^d} f\,\hf_\al.
\end{align*}

$\left.\right.$\vspace{-1.5pc}

\hfill $\Box$
\end{proof}

\begin{lem}\label{lem:xaLasobre2}
Let $p\in [1,\infty]$ and $\ord>0${\rm ,} then the operator
\begin{equation}\label{xaHaop}
|x|^{2\ord} H^{-\ord}f
\end{equation}
is bounded on $L^p(\RR^d)$.
\end{lem}

\begin{proof}
We will see that the kernel of \eqref{xaHaop} satisfies
\begin{equation}\label{inteny}
|x|^{2\ord}\int_{\RR^n} K_\ord(x,y)\,\hbox{d}y\,\leq\,C
\end{equation}
and
\begin{equation}\label{intenx}
\int_{\RR^n} |x|^{2\ord}K_\ord(x,y)\,\hbox{d}x\,\leq\,C,
\end{equation}
where the constant $C$ depends only on $d$ and $\ord$. Thus,
$|x|^{2\ord} H^{-\ord}f$ is bounded on $L^p(\RR^d)$ for all $p\in
[1,\infty]$ (Theorem~6.18 in \cite{4}).

We deal first with \eqref{inteny}. From \eqref{acotPhi} and
\eqref{laPhia} we see that there exists a constant $C$ such that
\eqref{inteny} is valid for all $|x|<2$.

Assume that $|x|<2|x-y|$.  By using eq.~\eqref{Kkercambio} we
obtain
\begin{align*}
K_\ord(x,y) &\leq
\frac{\hbox{e}^{-\frac{|x-y|^2}{8}}}{\Gamma(\ord)(4\pi)^{d/2}2^{\ord-1}}
\int_0^1 \zeta_\ord (s)\,
\hbox{e}^{-\frac{1}{8}\left(s|x+y|^2+\frac{1}{s}|x-y|^2\right)}\,
\frac{\hbox{d}s}{s}\\[.4pc]
&= \hbox{e}^{-\frac{|x-y|^2}{8}}
K_\ord\left(\frac{x}{\sqrt{2}},\frac{y}{\sqrt{2}}\right).
\end{align*}
Therefore, for some constant $C$,
\begin{align*}
\int_{\{|x|<2|x-y|\}} |x|^{2\ord}K_\ord(x,y)\,\hbox{d}y &\leq
\int_{\RR^d} |x-y|^{2\ord} \hbox{e}^{-\frac{|x-y|^2}{8}}
K_\ord\left(\frac{x}{\sqrt{2}},\frac{y}{\sqrt{2}}\right)\,\hbox{d}y\\[.4pc]
&\leq C\,\int_{\RR^d}
K_\ord\left(\frac{x}{\sqrt{2}},\frac{y}{\sqrt{2}}\right)\,\hbox{d}y,
\end{align*}
which is bounded by a constant independent of $x$.

It remains to consider integral \eqref{inteny} restricted to the
set $E_x=\{y\hbox{:}\ |x|>2|x-y|\}$ when $|x|>2$. Observe that in
this part, due to the identity $|x+y|^2=2|x|^2 -|x-y|^2+2|y|^2$,
we have
\begin{equation}\label{x_menor_xmasy}
|x|<|x+y|.
\end{equation}

As in the proof of Proposition~\ref{propo:KerHal}, we consider
$K_\ord = K_{\ord,0}+K_{\ord,1}$. Then, by using
\eqref{x_menor_xmasy},
\begin{equation*}
|x|^{2\ord} \int_{E_x} K_{\ord,1}(x,y) \hbox{d}y \le
C|x|^{2\ord}\hbox{e}^{-\frac18|x|^2}\int_{\RR^d}
\hbox{e}^{-\frac14|x-y|^2}\hbox{d}y \le C,
\end{equation*}
where in the last inequality we have used that for each positive
$b$, there exists a constant $C_b$ such that $|x|^b
\hbox{e}^{-|x|} \le C_b$.

In order to handle $K_{\ord,0}$, from \eqref{Z_0y1med} and
\eqref{x_menor_xmasy}, after some changes of variables we obtain
\begin{align*}
\int_{E_x} K_{\ord,0}(x,y) \hbox{d}y &\leq C \int_{E_x}
\int_0^{1/2} s^{-\frac{d}{2} +\ord}
\hbox{e}^{-\frac14\left(s|x|^2+ \frac1{s}|x-y|^2\right)}\frac{\hbox{d}s}{s} \hbox{d}y\\[.5pc]
&= C\, |x|^{d-2a}\int_{E_x} \int_0^{{|x|^2}/{2}} u^{-\frac{d}{2}
+\ord}
\hbox{e}^{-\frac14\left(u+\frac1{u}(|x||x-y|)^2\right)}\frac{\hbox{d}u}{u} \hbox{d}y\\[.5pc]
&= C\, |x|^{-2a}\int_0^{{|x|^2}/{2}} \int_0^{{|x|^2}/{2}}
u^{-\frac{d}{2} +\ord}
\hbox{e}^{-\frac14\left(u+\frac{r^2}{u}\right)}\frac{\hbox{d}u}{u}r^d
\frac{\hbox{d}r}{r}.
\end{align*}
The last double integral is bounded by
\begin{align*}
\int_0^{\infty} \int_0^{\infty} u^{-\frac{d}{2} +\ord}
\hbox{e}^{-\frac14\left(u+\frac{r^2}{u}\right)}\frac{\hbox{d}u}{u}\,r^d\,
\frac{\hbox{d}r}{r} &= \int_0^{\infty} u^{-\frac{d}{2} +\ord}
\hbox{e}^{-\frac{u}{4}} \int_0^{\infty}
\hbox{e}^{-\frac{r^2}{4u}}\, r^d\, \frac{\hbox{d}r}{r}\,
\frac{\hbox{d}u}{u}\\[.5pc]
&= \int_0^{\infty} \hbox{e}^{-\frac{u}{4}}\,
u^{\ord}\,\frac{\hbox{d}u}{u} \int_0^{\infty}
\hbox{e}^{-\frac{r^2}{4}}\, r^d\, \frac{\hbox{d}r}{r},
\end{align*}
and both integrals clearly converge.\pagebreak

Finally we shall prove \eqref{intenx}. We split
\begin{equation}\label{ultima}
\int_{\RR^d} |x|^{2\ord} K_\ord(x,y)\,\hbox{d}x =
\left(\int_{\{|y|>2|x-y|\}} + \int_{\{|y|<2|x-y|\}}\right)
|x|^{2\ord} K_\ord(x,y)\,\hbox{d}x.
\end{equation}
Since $|x|\leq \frac{3}{2}|y|$ when $|y|>2|x-y|$ and the kernel
$K_\ord$ is symmetric, the first integral of \eqref{ultima} is
less than or equal to
\begin{equation*}
\left(\frac{3}{2}\,|y|\right)^{2\ord}\int_{\{|y|>2|x-y|\}}\!\!
K_\ord(y,x)\,\hbox{d}x,
\end{equation*}
which can be bounded, in the same way as \eqref{inteny}, by a
constant depending only on $d$ and $\ord$. For the second term of
\eqref{ultima}, we have
\begin{equation*}
\int_{\{|y|<2|x-y|<4\}} |x|^{2\ord} K_\ord(x,y)\,\hbox{d}x\,
\leq\, C\int_{\RR^d} K_\ord(x,y)\,\hbox{d}x.
\end{equation*}
On the other hand, by estimate \eqref{acotPhi} and the expression
\eqref{laPhia},
\begin{align*}
\int_{\{|y|<2|x-y|,\ 4<|x-y|\}} |x|^{2\ord} K_\ord(x,y)\,\hbox{d}x
&\leq
C\,\int_{\{4<|x-y|\}} |x-y|^{2\ord} \hbox{e}^{-\frac{|x-y|^2}{8}} \,\hbox{d}x\\[.4pc]
&\leq C\,\int_{\{4<|x|\}} |x|^{2\ord} \hbox{e}^{-\frac{|x|^2}{8}} \,\hbox{d}x\\[.4pc]
&\leq C.
\end{align*}
Then, we have proved \eqref{intenx}.\hfill $\Box$
\end{proof}

\section{Potential spaces}

\begin{defini}\label{defi:esppot}$\left.\right.$\vspace{.5pc}

\noindent {\rm Given $p\in[1,\infty)$ and $\ord> 0$, we define the
space
\begin{equation*}
\LL_\ord^p = H^{-\ord/2}(L^p(\RR^d))\text{,}
\end{equation*}
with a norm given by
\begin{equation*}
\|f\|_{\LL_\ord^p} = \|g\|_p,
\end{equation*}
where $g$ is such that $H^{-\ord/2}g = f$.}
\end{defini}

\begin{rem}{\rm
The space $\LL^p_\ord$ is well-defined for $p\in[1,\infty)$ and
$\ord> 0$, since $H^{-\ord/2}$ is one-to-one. As
$\D=H^{-\ord/2}(\D)$ then $\D$ is a dense space of $\LL^p_\ord$.}
\end{rem}

Due to the boundedness of the operator $H^{-\ord}$, the space
$\LL^p_\ord$ is a subspace of $L^p(\RR^d)$. Moreover, we have the
following theorem.\pagebreak

\begin{theor}[\!]\label{teo:isom_orden}
Let $0< \ord<b${\rm ,} then
\begin{enumerate}
\renewcommand\labelenumi{\rm (\roman{enumi})}
\leftskip .15pc
\item $\LL^p_b \subset \LL^p_\ord \subset L^p$ and the inclusions
are continuous.

\item $\LL^p_\ord$ and $\LL^p_b$ are isometrically isomorphic.
\end{enumerate}
\end{theor}

\begin{proof}
Since $H^{-b/2}= H^{-\ord/2}\circ H^{-\ga}$, where
$\ga=(b-\ord)/2$, we see that (i) easily follows from the
definition of $\LL^p_\ord$ and $\LL^p_b$, and the boundedness of
$H^{-\ga}$. From Definition~\ref{defi:esppot} and the fact that
$H^{-\ga}$ is one-to-one, it is easy to verify that
$H^{-\ga}\hbox{:}\ \LL^p_\ord \mapsto \LL^p_b$ is an isometric
isomorphism, and this gives (ii).\hfill $\Box$
\end{proof}

The members of $\LL^p_\ord$ have a special decay at infinity, as
the following proposition shows.

\begin{propo}\label{propo:xaf_en_La}$\left.\right.$\vspace{.5pc}

\noindent If $p\in [1,\infty)${\rm ,} $\ord>0$ and $f\in
\LL^p_\ord${\rm ,} then $|x|^\ord f(x)$ belongs to $L^p(\RR^d)$.
\end{propo}

\begin{proof}
By Definition~\ref{defi:esppot}, we have $f= H^{-\ord/2}g$ with $g
\in L^p(\RR^d)$. Then it is enough to apply
Lemma~\ref{lem:xaLasobre2}.\hfill $\Box$
\end{proof}

We remind the classical Sobolev spaces $L^p_\ord$, with
$p$ in the range $1\leq p<\infty$ and $\ord>0$, defined by
\begin{equation}\label{sob:clasicos}
L^p_\ord = (I-\Delta)^{-\ord/2}(L^p(\RR^d)).
\end{equation}
In this case, the norm for $f\in L^p_\ord$ is given by
\begin{equation*}
\|f\|_{L^p_\ord} = \|g\|_p,
\end{equation*}
where the function $g\in L^p$ is such that
$(I-\Delta)^{-\ord/2}g=f$ (see \cite{8}). In the following theorem
we describe the relation between the spaces $L^p_\ord$ and
${\LL}^p_\ord$.

\begin{theor}[\!]\label{teo:HSvsS}

Let $\ord>0$ and $p\in(1,\infty)${\rm ,} then
\begin{enumerate}
\renewcommand\labelenumi{\rm (\roman{enumi})}
\leftskip .4pc
\item \label{contenido}$\LL^p_\ord \subset L^p_\ord$.

\item \label{noigual}$\LL^p_\ord \neq L^p_\ord$.

\item \label{Sisuppcomp} If $f\in L^p_\ord$ and has compact
support{\rm ,} then $f$ belongs to $\LL^p_\ord$.
\end{enumerate}
\end{theor}

\begin{proof}
To see \eqref{contenido} we follow the argument in \cite{14}.
Namely, the symbol of
\begin{equation*}
(I-\Delta)^{\ord/2}H^{-\ord/2}
\end{equation*}
belongs to the class $S^0_{1,0}$ and so it defines a bounded
operator on $L^p(\RR^d)$ (see \cite{9}). Let $f\in\LL^p_\ord$,
then $h=[(I-\Delta)^{\ord/2}H^{-\ord/2}](H^{\ord/2}f)$ is a
function of $L^p(\RR^d)$ with $(I-\Delta)^{-\ord/2}h=f$. Hence $f$
belongs to $L^p_\ord$.

In order to see (ii) let
\begin{equation*}
g(x)=\frac{1}{(1+|x|)^{1/p+\ord}}\quad \hbox{and}\quad
f=(I-\Delta)^{-\ord/2}g.
\end{equation*}
Since $g\in L^p(\RR^d)$, then $f\in L^p_\ord$. We will see that
$f$ is not in $\LL_\ord^p$. From
Proposition~\ref{propo:xaf_en_La}, if $f$ were in $\LL_\ord^p$ we
would have $|x|^\ord f\in L^p(\RR^d)$. However, if
$\mathcal{G}_{\ord}$ is the kernel of $(I-\Delta)^{-\ord/2}$,
\begin{equation*}
\mathcal{G}_\ord(x)=\frac{1}{(4\pi)^{\ord/2}\Gamma(\ord/2)}\int_0^\infty
\hbox{e}^{-\frac{\pi\,|x|^2}{t}-\frac{t}{4\pi}}\,t^{\frac{-d+\ord}{2}}\frac{\hbox{d}t}{t},
\end{equation*}
then
\begin{align*}
f(x) &= (I-\Delta)^{-\ord/2} g(x) = \int_{\RR^d} \mathcal{G}_\ord(y) g(x-y)\,\hbox{d}y\\[.3pc]
&\ge \int_{\{|y|<1\}}\! \mathcal{G}_\ord(y)\, g(x-y)\,\hbox{d}y
\ge (2+|x|)^{-1/p-\ord}\!\int_{\{|y|<1\}}\!
\mathcal{G}_\ord(y)\,\hbox{d}y.
\end{align*}
Thus, $|x|^\ord f$ is not in $L^p(\RR^d)$.

Finally, (iii) is a direct consequence of the following
inequality,
\begin{equation*}
\int_{\RR^d} |H^\ord(\psi f)|^p \leq C(\psi) \int_{\RR^d}
|(-\Delta + I)^\ord f|^p,
\end{equation*}
where $\psi\in \CC$ and $C(\psi)$ is a constant depending on
$\psi$. This inequality was proved in \cite{14} for $p=2$ and the
same proof works for $p$ in the range $1<p<\infty$.\hfill $\Box$
\end{proof}

\begin{theor}[\!]\label{teo:Wkp}
Let $k\in \mathbb{N}$ and $p\in (1,\infty)${\rm ,} then
\begin{equation*}
W^{k,p}=\LL^p_k
\end{equation*}
and the norms $\|\cdot\|_{W^{k,p}}$ and $\|\cdot\|_{\LL^p_k}$ are
equivalent.
\end{theor}

We first present some technical results that we shall need for the
proof of this theorem.

\begin{lem}\label{lem:conmut}
Let $b\in\RR${\rm ,} then for all $f$ in $\D${\rm ,} we have
\begin{align}
A_jH^{b}f &= (H+2)^{b}A_jf,\quad 1\leq j \leq d,\label{conmut+}\\[.3pc]
A_jH^{b}f &= (H-2)^{b}A_jf,\quad -d\leq j \leq-1,\label{conmut-}\\[.3pc]
H^{b}A_jf &= A_j(H-2)^{b}f,\quad 1\leq j \leq d,\label{conmut+p}
\end{align}
and
\begin{equation}\label{conmut-p}
H^{b}A_jf = A_j(H+2)^{b}f\text{,}\ \ \ -d\leq j \leq -1,
\end{equation}
where $H^b\hf_\al=(2|\al|+d)^b\hf_\al$ and
$(H+2)^b\hf_\al=(2|\al|+d+2)^b\hf_\al${\rm ,} for all
$\al\in\NN^d${\rm ,} and
$(H-2)^b\hf_\al=(2|\al|+d-2)^b\hf_\al${\rm ,} for all
$\al\in\NN^d$ with $|\al|\ge 1$.
\end{lem}

\begin{proof}
Let $1\leq j \leq d$ and $\al\in \NN^d$, then
\begin{align*}
A_jH^b\hf_{\al} &= (2|\al|+d)^b\, A_j\hf_{\al} =
\sqrt{2\al_j}\,(2|\al|+d)^b\, \hf_{\al-e_j}\\[.4pc]
&= \sqrt{2\al_j}\,(2(|\al|-1)+d+2)^b\, \hf_{\al-e_j} =
\sqrt{2\al_j}\,(H+2)^b \hf_{\al-e_j}\\[.4pc]
&= (H+2)^b A_j \hf_{\al}\,,
\end{align*}
and this gives \eqref{conmut+} by linearity. In the same way we
obtain \eqref{conmut-} and \eqref{conmut-p}. We are assuming in
\eqref{conmut+p} that $f$ is a linear combination of Hermite
functions with order $|\al|\ge 1$.\hfill $\Box$
\end{proof}

The Hermite--Riesz transforms associated to $H$ are defined as
\begin{equation*}
R_j = A_j H^{-1/2},\quad 1\leq |j| \leq d
\end{equation*}
and the Hermite--Riesz transform vector
\begin{equation*}
R = (R_{-d}, R_{-d+1},\dots,R_{-1},R_{1},\dots,R_{d-1},R_{d}).
\end{equation*}

These operators were introduced by Thangavelu in \cite{12} (see
also \cite{10}). He proved that they are bounded on $L^p(\RR^d)$
for $1<p<\infty$, and of weak type $(1,1)$. Given $m\in
\mathbb{N}$ the Hermite--Riesz of order $m$ is defined as
\begin{equation}\label{Rsup}
R_{j_1, j_2,\dots, j_m} = A_{j_1}A_{j_2}\cdots A_{j_m} H^{-m/2},
\end{equation}
where $1\leq |j_n|\leq d$, for every $1\leq n\leq m$.

For further references, we enounce the following crystallized
theorem of known facts about the Hermite--Riesz transforms.

\begin{theor}[\!]\label{teo:cristal}$\left.\right.$
\begin{enumerate}
\renewcommand\labelenumi{\rm (\alph{enumi})}
\leftskip .1pc
\item \label{Rpseudo} The Hermite--Riesz transforms $R_j${\rm ,}
$1\leq |j|\leq d${\rm ,} are pseudo-differential operators whose
symbols belong to $S^0_{1,0}$. In particular{\rm ,} they are
bounded in the classical Sobolev spaces $\eqref{sob:clasicos}$.

\item \label{acotRvect} Let $m\in \mathbb{N}$ and $p\in
(1,\infty)$. Then there exists a constant $C_{p,m}$ not depending
on the dimension $d${\rm ,} such that
\begin{equation*}
\hskip -1.25pc \left\| \left(\sum_{1\leq |j_1|,\dots,|j_m|\leq d}
|R_{j_1,\dots,j_m}f|^2 \right)^{1/2}\right\|_p\ \leq C_{p,m}\,
\|f\|_p.
\end{equation*}
\end{enumerate}
\end{theor}

\begin{proof}
For the proof of (a) see \cite{13}. For (b) see \cite{5} (see also
\cite{10} and \cite{12} for the case $m=1$).\hfill $\Box$
\end{proof}

Now we have the following proposition.

\begin{propo}\label{propo:Reves}$\left.\right.$\vspace{.5pc}

\noindent For $f$ and $g$ in $\D${\rm ,} we have
\begin{equation*}
\int_{\RR^d} fg = 2 \int_{\RR^d} \sum_{1\leq |j|\leq d} R_{j}f\,
R_{j}g.
\end{equation*}
Let $p\in(1,\infty)$. Then there exists a constant $C$ such that
for all $f$ in $\D${\rm ,} we have
\begin{equation}\label{Reves}
\|f\|_p \leq C\,\left\| \left(\sum_{1\leq |j|\leq d} |R_{j}f|^2
\right)^{1/2}\right\|_p.
\end{equation}
\end{propo}

\begin{proof}
For $1\leq |j|\leq d$, we have
\begin{equation*}
R_j^*R_j=H^{-1/2}A_j^*A_jH^{-1/2}=H^{-1/2}A_{-j}A_jH^{-1/2}.
\end{equation*}
Then by formula \eqref{Hfact},
\begin{equation*}
\sum_{1\leq |j|\leq d} R_j^*R_j = H^{-1/2}\left(\sum_{1\leq
|j|\leq d}A_{-j}A_j\right)H^{-1/2} = 2I.
\end{equation*}

Therefore, if $f$ and $g$ are in $\D$,
\begin{align*}
\sum_{1\leq |j|\leq d} \int_{\RR^d} R_j f \,R_j g &= \sum_{1\leq
|j|\leq d} \int_{\RR^d} R_j^* R_j f g = \int_{\RR^d}
\left(\sum_{1\leq |j|\leq d}R_j^*R_jf\right)\,g\\[.3pc]
&= 2\,\int_{\RR^d} fg.
\end{align*}

In order to prove \eqref{Reves}, by H\"older's inequality, we get
\begin{align*}
\hskip -4pc \|f\|_p &= \sup_{\{g\in \D\hbox{:}\ \|g\|_{p'}=1\}}
\int_{\RR^d} fg = \frac{1}{2}\sup_{\{g\in \D\hbox{:}\
\|g\|_{p'}=1\}}\ \sum_{1\leq |j|\leq
d} \int_{\RR^d} R_j(f)R_j(g)\\[.5pc]
\hskip -4pc &\leq \frac{1}{2}\sup_{\{g\in \D\hbox{:}\
\|g\|_{p'}=1\}}\ \left\| \left(\sum_{1\leq |j| \leq d}
|R_{j}(f)|^2 \right)^{1/2}\right\|_p\, \left\| \left(\sum_{1\leq
|j| \leq d}
|R_{j}(g)|^2 \right)^{1/2}\right\|_{p'}\\[.5pc]
\hskip -4pc &\leq C \left\| \left(\sum_{1\leq |j| \leq d}
|R_{j}(f)|^2 \right)^{1/2}\right\|_p,
\end{align*}
where in the last inequality we have used
Theorem~\ref{teo:cristal}(b).\hfill $\Box$
\end{proof}

\setcounter{poft}{3}
\begin{poft}{\rm
Since $\D$ is dense in both spaces it is enough to show the
equivalence of the norm for functions in $\D$. Let $f\in \D$, and
let $f=H^{-k/2}g$. Then by Theorem~\ref{teo:cristal}(b) and
Theorem~\ref{teo:fcoef}, we obtain
\begin{align*}
\|f\|_{W^{k,p}} &= \sum_{1\leq |j_1|,\dots,|j_k| \leq d ,\ 1\leq
m\leq k} \|R_{j_1,\dots,j_m}H^{-(k-m)/2}g\|_p + \|H^{-k/2}g\|_p\\[.3pc]
&\leq C_{k,d}\,\|g\|_p = C_{k,d}\,\|f\|_{\LL^p_k}.
\end{align*}

To prove the converse inequality, we first consider the case
$k=1$. By Proposition~\ref{propo:Reves}, we get
\begin{equation*}
\|f\|_{\LL^p_1} = \|H^{1/2}f\|_p \leq C \sum_{1\leq |j|\leq d}
\|A_jf\|_p \leq C\, \|f\|_{W^{1,p}}.
\end{equation*}

Now we shall use an inductive argument. Suppose we have
\begin{equation*}
\|f\|_{\LL^p_m} = \|H^{m/2}f\|_p \leq C_{m} \|f\|_{W^{m,p}},
\end{equation*}
for $f\in \D$, with $m<k$. Since for some constants $c_1$,
$c_2$,\dots,$c_{k-1}$,
\begin{equation*}\label{identH}
\sum_{1\leq |j_1|,\dots,|j_k|\leq d} A^*_{j_k}\cdots
A^*_{j_1}\,A_{j_1}\cdots A_{j_k} = 2^k H^k + \sum_{m=1}^{k-1}
c_m\, H^m,
\end{equation*}
and, since $H$ is autoadjoint, for all $f, g\in \D$,
\begin{align*}
\|f\|_{\LL^p_k} &= \|H^{k/2}f\|_p = \sup_{\{g\in \D\hbox{:}\
\|g\|_{p'}= 1\}} \int_{\RR^d} (H^{k/2}f)\,g\\[.3pc]
&= \sup_{\{g\in \D\hbox{:}\ \|g\|_{p'}=1\}}  \int_{\RR^d}
(H^{k}f)\,(H^{-k/2}g).
\end{align*}

Now by using formula \eqref{identH} and the definition
\eqref{Rsup} of the Hermite--Riesz transform of higher order, we
have
\begin{align*}
&2^k\int_{\RR^d} H^{k}f\,H^{-k/2}g\\[.3pc]
&= \int_{\RR^d} \left(\sum_{1\leq |j_1|,\dots,|j_k|\leq d}
A^*_{j_k}\cdots A^*_{j_1}A_{j_1}\cdots A_{j_k} - \sum_{m=1}^{k-1}
c_m\, H^m\right) f\,H^{-k/2}g\\[.3pc]
&=\! \sum_{1\leq |j_1|,\dots,|j_k|\leq d}
\int_{\RR^d}\!A_{j_1}\cdots A_{j_k}f R_{j_1\dots,j_k}g -
\sum_{m=1}^{k-1} c_m\, \int_{\RR^d} H^{m/2}\!f\,
H^{-\frac{(k-m)}{2}}\!g.
\end{align*}
Thus, by H\"{o}lder's inequality the last expression is bounded by
\begin{align}\label{sumadenorm}
&\sum_{1\leq |j_1|,\dots,|j_k|\leq d} \|A_{j_1}\cdots A_{j_k}f\|_p
\|R_{j_1\dots,j_k}g\|_{p'}\nonumber\\[.3pc]
&\quad\, + \sum_{m=1}^{k-1} |c_m| \,\|H^{m/2}f\|_p
\|H^{-\frac{(k-m)}{2}}g\|_{p'}.
\end{align}
From Theorem~\ref{teo:cristal}(b), Theorem~\ref{teo:fcoef}, the
induction hypotheses and Definition~\ref{defi:sob}, there exists a
constant $C$ such that \eqref{sumadenorm} is bounded by
\begin{equation*}
\left(C\sum_{m=1}^{k-1} |c_m|\right) \|f\|_{W^{k,p}}.
\end{equation*}}
\end{poft}

\section{Boundedness of some operators on
$\pmb{\LL^p_\ord}$}\label{sec:boud_op_LL}

\begin{theor}[\!]\label{teo:acotAj}
Let $p\in(1,\infty)${\rm ,} $\ord> 1$ and $1\leq |j| \leq d$. Then
$A_j$ is bounded from $\LL^p_\ord$ into $\LL^p_{\ord-1}$.
\end{theor}

\begin{proof}
Let $1\leq j\leq d$ (the case $-d\leq j\leq -1$ is similar). If
$f\in \D$, by Lemma~\ref{lem:conmut},
\begin{equation}\label{paraAj}
A_jf= H^{-(\ord-1)/2}\left(\frac{H}{H+2}\right)^{(\ord-1)/2}\!
R_j\, H^{\ord/2}f.
\end{equation}
As the function $m(\la)= ({\la}/({\la+2}))^{(a-1)/2}$ satisfies
the hypotheses of Theorem~4.2.1 in \cite{12}, the operator
$({H}/({H+2}))^{(\ord-1)/2}$ is bounded on $L^p(\RR^d)$. Hence, by
Theorem~\ref{teo:cristal}(a) and Theorem~\ref{teo:fcoef}, the
operator
\begin{equation}\label{opHR}
H^{-(\ord-1)/2}\left(\frac{H}{H+2}\right)^{(\ord-1)/2}\! R_j
\end{equation}
is bounded on $L^p(\RR^d)$. If $f\in \LL^p_\ord$, then $A_jf=R_j
H^{-(\ord-1)/2}H^{\ord/2}f$. Since operators \eqref{opHR} and $R_j
H^{-(\ord-1)/2}$ coincide in $\D$, both are bounded on
$L^p(\RR^d)$ and $\D$ is dense in $\LL^p_\ord$, formula
\eqref{paraAj} also works for all $f\in\LL^p_\ord$.

Therefore, if $f\in\LL^p_\ord$, the function
$h=({H}/{(H+2)})^{(\ord-1)/2} R_j\, H^{\ord/2} f$ belongs to
$L^p(\RR^d)$ and by \eqref{paraAj}, we have
\begin{equation*}
\|A_jf\|_{\LL^p_{\ord-1}}=\|h\|_p\leq C\|H^{\ord/2} f\|_p=
C\|f\|_{\LL^p_\ord}.
\end{equation*}

$\left.\right.$\vspace{-1.5pc}

\hfill $\Box$
\end{proof}

\begin{theor}[\!]
If $p\in (1,\infty)${\rm ,} $\ord>0$ and $1\leq |j|\leq d${\rm ,}
then the operators $R_j$ are bounded on $\LL_\ord^p$.
\end{theor}

\begin{proof}
By Theorem~\ref{teo:isom_orden}, $H^{-1/2}$ is bounded from
$\LL^p_\ord$ into $\LL^p_{\ord+1}$. Hence,
Theorem~\ref{teo:acotAj} gives the desired result.\hfill $\Box$
\end{proof}

\begin{rem}{\rm
As the Hermite--Riesz transforms are pseudo-differential operators
with symbols in the class $S^0_{1,0}$ (see
Theorem~\ref{teo:cristal}(a)), they map the classical Sobolev
spaces $L^p_\ord$ into themselves for all $p\in(1,\infty)$ and
$\ord>0$. However, the classical Riesz transforms are not bounded
on $\LL^p_\ord$ for any $\ord>({1}/{p'})$ as the following
proposition shows.}
\end{rem}

\begin{propo}$\left.\right.$\vspace{.5pc}

\noindent Let $p$ be in the range $1<p<\infty.$ The Hilbert
transform on $\RR$ is not bounded in $\LL^p_\ord$ for any
$\ord>({1}/{p'})$.
\end{propo}

\begin{proof}
Let $\mathcal{H}$ be  the Hilbert transform on the line, that is
\begin{equation*}
\mathcal{H}f(x) = \lim_{\eps\rightarrow 0} \int_{|y|>\eps}\!
\frac{f(x-y)}{y}\,\hbox{d}y.
\end{equation*}
Consider the function $f$ in $\LL^p_\ord$ given by
\begin{equation*}
f(x)= \exp(-|x|^2).
\end{equation*}
Given $\eps>0$ and $x>2$, we have
\begin{equation}\label{trunc}
\int_{|y|>\eps}\! \frac{f(x-y)}{y}\,\hbox{d}y =
\left(\int_{\eps<|y|<1} + \int_{-\infty}^{-1} +
\int_1^\infty\right)\! \frac{f(x-y)}{y}\,\hbox{d}y.
\end{equation}
By the mean value theorem, the first integral of the last
expression can be written as
\begin{align*}
\int_{\eps<|y|<1}\! \frac{\hbox{e}^{-|x-y|^2}}{y}\,\hbox{d}y &= -
\int_{\eps<|y|<1}\!
\frac{\hbox{e}^{-|x|^2}-\hbox{e}^{-|x-y|^2}}{x-(x-y)}\,\hbox{d}y\\[.4pc]
&= -2 \int_{\eps<|y|<1}\!
\theta(x,y)\,\hbox{e}^{-\theta(x,y)^2}\,\hbox{d}y,
\end{align*}
where $x-1<x-y<\theta(x,y)<x$ if $\eps<y<1$, and
$x<\theta(x,y)<x-y<x+1$ if $-1<y<-\eps$. Thus
\begin{align}\label{1int}
\left|\int_{\eps<|y|<1}\!
\frac{\hbox{e}^{-|x-y|^2}}{y}\,\hbox{d}y\right| &< 2
\left|\int_{\eps<|y|<1}\!
\theta(x,y)\,\hbox{e}^{-\theta(x,y)^2}\,\hbox{d}y\right|\nonumber\\[.4pc]
&\leq\ 4(x+1)\, \hbox{e}^{-|x-1|^2}.
\end{align}
For the second integral of \eqref{trunc},
\begin{equation}\label{2int}
\left|\int_{-\infty}^{-1}\!
\frac{\hbox{e}^{-|x-y|^2}}{y}\,\hbox{d}y\right|\ \leq\
\hbox{e}^{-|x|^2} \int_{-\infty}^{-1}\!
\frac{\hbox{e}^{-|y|^2}}{|y|}\,\hbox{d}y.
\end{equation}
Finally, since $x>2$,
\begin{equation}\label{3int}
\left|\int^{\infty}_{1}\!
\frac{\hbox{e}^{-|x-y|^2}}{y}\,\hbox{d}y\right|\ \ge\ \frac{1}{x}
\int^{x}_{1}\! \hbox{e}^{-|x-y|^2}\,\hbox{d}y\ \ge\ \frac{1}{x}
\int^{1}_{0}\! \hbox{e}^{-|u|^2}\,\hbox{d}u.
\end{equation}
Therefore, from \eqref{1int}--\eqref{3int}, there exist constants
$C$ and $M$ independent of $\eps$, such that
\begin{equation*}
\left|\int_{|y|>\eps}\! \frac{f(x-y)}{y}\,\hbox{d}y\right| \ge
\frac{C}{x}\quad \text{for all $x>M$.}
\end{equation*}

Thus, for any $\ord>({1}/{p'})$, $|x|^{\ord}\, \mathcal{H}f$ is
not in $L^p(\RR)$, and by Proposition~\ref{propo:xaf_en_La},
$\mathcal{H}f$ is not in $\LL^p_\ord$.\hfill $\Box$
\end{proof}

\section{Poincar\'e inequalities}

This section is devoted to some Poincar\'e-type inqualities.

\begin{rem}\label{rem:acotH_int}{\rm
Observe that by formula \eqref{laPhia}, if $\ord<d$ and we take a
positive function $f$, then
\begin{equation*}
H^{-\ord/2}f(x)\leq C\,\int_{\RR^d}
\frac{f(y)}{|x-y|^{d-\ord}}\hbox{d}y.
\end{equation*}
Therefore, in the case $\ord<d$ the operator $H^{-\ord/2}$
inherits the boundedness properties of the fractional integral. In
particular, $H^{-\ord/2}$ is bounded from $L^p(\RR^d)$ into
$L^q(\RR^d)$ for $\frac{1}{q}=\frac{1}{p}-\frac{\ord}{d}$, where
$1< p < q < \infty$.}
\end{rem}

Next theorem gives the behavior of $H^{-\ord/2}$ on $L^p(\RR^d)$,
for $p\ge 1$. Inequality \eqref{inteny} allows us to obtain some
boundedness of the operator $H^{-\ord/2}$ that has a different
flavor from the boundedness of the classical fractional integral.

\begin{theor}[\!]\label{teo:acotHal}
Let $a,d$ such that $0<a<d${\rm ,} then
\begin{enumerate}
\renewcommand\labelenumi{\rm (\roman{enumi})}
\leftskip .4pc
\item There \label{acotL1Lq} exists a constant $C${\rm ,} such that
\begin{equation*}
\hskip -1.25pc \|H^{-\ord/2}f\|_{q} \leq C \|f\|_1,
\end{equation*}
for all $f\in L^1(\RR^d)$ if and only if $1\leq q<{d}/({d-\ord})$.

\item \label{acotLpLinf} There exists a constant $C${\rm ,} such that
\begin{equation*}
\hskip -1.25pc \|H^{-\ord/2}f\|_{\infty} \leq C \|f\|_{p},
\end{equation*}
for all $f$ in $L^p(\RR^d)$ if and only if $p>\frac{d}{\ord}$.

\item \label{item:acotLinfLq} There exists a constant $C${\rm ,} such that
\begin{equation*}
\hskip -1.25pc \|H^{-\ord/2}f\|_{q} \leq C \|f\|_{\infty},
\end{equation*}
for all $f\in L^\infty(\RR^d)$ if and only if $q>({d}/{\ord})$.

\item \label{item:acotLpL1} There exists a constant $C${\rm ,} such that
\begin{equation*}
\hskip -1.25pc \|H^{-\ord/2}f\|_{1} \leq C \|f\|_{p},
\end{equation*}
for all $f\in L^p(\RR^d)$ if and only if $1\leq p <
{d}/({d-\ord})$.

\item \label{item:acotHal} If $1< p< \infty${\rm ,} $1< q< \infty$ and
$\frac{1}{p}-\frac{\ord}{d}\leq
\frac{1}{q}<\frac{1}{p}+\frac{\ord}{d}${\rm ,} then there exists a
constant $C$, such that
\begin{equation*}\label{acotHalLp}
\hskip -1.25pc \|H^{-\ord/2}f\|_{q} \leq C \|f\|_p,
\end{equation*}
for all $f\in L^p(\RR^d)$.
\end{enumerate}
\end{theor}

\begin{proof}
By Minkowski's integral inequality,
\begin{align*}
\int_{\RR^d} (H^{-\ord/2} f)^q &= \int_{\RR^d}
\left(\int_{\RR^d}K_{\ord/2}(x,y) f(y)\,\hbox{d}y\right)^q \,\hbox{d}x\\[.3pc]
&\leq \left(\int_{\RR^d} f(y)\, \left(\int_{\RR^d}
K_{\ord/2}(x,y)^q \,\hbox{d}x\right)^{1/q} \,\hbox{d}y\right)^q,
\end{align*}
for all $f\in L^1(\RR^d)$. From inequalities \eqref{acotPhi} and
\eqref{laPhia}, we get
\begin{align*}
\int_{\RR^d} K_{\ord/2}(x,y)^q \,\hbox{d}x\leq C\,
\left(\int_{\{|x|<2\}} \frac{\hbox{d}x}{|x|^{q(d-\ord)}} +
\int_{\{|x|>2\}} \hbox{e}^{-q\frac{|x|^2}{4}}\,\hbox{d}x\right),
\end{align*}
and this integral is finite if $1\leq q<{d}/({d-\ord})$, proving
(i). Conversely, by \eqref{Z_0y1med} and a change of variables, if
$|x-y|<1$,
\begin{align}\label{abajo}
K_{\ord/2}(x,y) &\ge
\frac{\hbox{e}^{-|x+y|^2}}{C_1\Gamma(\ord/2)(4\pi)^{d/2}
2^{\ord/2-1}}\int_0^{1/2} s^{-\frac{d-\ord}{2}}\,
\hbox{e}^{-\frac{1}{4s}|x-y|^2}\, \frac{\hbox{d}s}{s}\nonumber\\[.3pc]
&\ge \frac{\int_0^{\frac{1}{2|x-y|^2}} s^{-\frac{d-\ord}{2}}
\hbox{e}^{-\frac{1}{4s}} \frac{{\rm
d}s}{s}}{C_1\Gamma(\ord/2)(4\pi)^{d/2}2^{\frac{\ord}{2}-1}}\
\frac{\hbox{e}^{-|x+y|^2}}{|x-y|^{d-\ord}}\nonumber\\[.3pc]
&\ge \frac{\int_0^{\frac{1}{2}} s^{-\frac{d-\ord}{2}}
\hbox{e}^{-\frac{1}{4s}} \frac{{\rm
d}s}{s}}{C_1\Gamma(\ord/2)(4\pi)^{d/2}2^{\frac{\ord}{2}-1}}\
\frac{\hbox{e}^{-|x+y|^2}}{|x-y|^{d-\ord}}.
\end{align}
Now, let $f_n$, $n\ge 0$, be an approximation to the identity.
Suppose that inequality (i) holds for all $f\in L^1(\RR^d)$, then
by inequality \eqref{abajo} there exists a constant $C$, such that
\begin{equation*}
C\ge \int_{\RR^d} \left(\int_{\RR^d}
\frac{\chi_{\{|x-y|<1\}}(y)\hbox{e}^{-|x+y|^2}}{|x-y|^{d-\ord}}
f_n(y)\,\hbox{d}y\right)^q \,\hbox{d}x,
\end{equation*}
for all $n\ge 0$, so that
\begin{equation*}
C\ge \int_{\{|x|<1\}}
\frac{\hbox{e}^{-q|x|^2}}{|x|^{q(d-\ord)}}\,\hbox{d}x,
\end{equation*}
but this is false when $q(d-\ord)\ge d$.\pagebreak

To obtain inequality (ii), by H\"older's inequality,
\begin{equation*}
\left|\int_{\RR^d} K_{\ord/2}(x,y)\,f(y)\,\hbox{d}y\right| \leq
\|f\|_p \left(\int_{\RR^d}
K_{\ord/2}(x,y)^{p'}\,\hbox{d}y\right)^{1/p'},
\end{equation*}
and in the same manner as we dealt with (i), the last integral is
finite when $p>({d}/{\ord})$. To see that $p>({d}/{\ord})$ is
necessary, let $\eps>0$ and
\begin{equation*}
f(x)= \begin{cases}
|x|^{-\ord}\left(\log\frac{1}{|x|}\right)^{-(\ord/d)(1+\eps)}
&\text{if $|x|\leq 1/2$,}\\[.3pc]
0 &\text{if $|x|>1/2$.}\end{cases}
\end{equation*}
Then $f\in L^p(\RR^d)$ for all $p\leq ({d}/{\ord})$. However,
$H^{-\ord}f$ is essentially unbounded as by estimate
\eqref{abajo},
\begin{equation*}
H^{-\ord/2}f(0)\ge C\int_{|x|\leq 1/2}
|x|^{-d}\left(\log\frac{1}{|x|}\right)^{-(\ord/d)(1+\eps)} =
\infty \text{,}
\end{equation*}
when $\eps$ is small enough.

To see (iii), let $f$ in $L^\infty(\RR^d)$, then
\begin{align*}
\int_{\RR^d} |H^{-\ord/2} f|^q & \leq \|f\|_\infty^q \int_{\RR^d}
\left(\int_{\RR^d}K_{\ord/2}(x,y)\,\hbox{d}y\right)^q \,\hbox{d}x\\[.3pc]
&\leq \|f\|_\infty^q \left[\int_{\{|x|\leq 1\}}\!\! + \int_{\{|x|>
1\}} \left(\int_{\RR^d}K_{\ord/2}(x,y)\,\hbox{d}y\right)^q
\,\hbox{d}x\right].
\end{align*}
The first integral of the last expression is finite due to
estimate \eqref{acotPhi}. By inequality \eqref{inteny}, the second
integral is
\begin{align*}
\int_{\{|x|> 1\}}
\left(\int_{\RR^d}K_{\ord/2}(x,y)\,\hbox{d}y\right)^q \,\hbox{d}x
&\leq C^q \int_{\{|x|> 1\}} \frac{\hbox{d}x}{|x|^{q\ord}},
\end{align*}
which is finite when $q>({d}/{\ord})$.

To see the converse, let $f\equiv 1$. Then,
\begin{align}\label{infinita}
\int_{\RR^d} |H^{-\ord/2} f|^q &= \int_{\RR^d}
\left(\int_{\RR^d}K_{\ord/2}(x,y)\,\hbox{d}y\right)^q \,\hbox{d}x,
\end{align}
and we will see that there exists a constant $c$ such that
\begin{align}\label{porabajo}
\int_{\RR^d}K_{\ord/2}(x,y)\,\hbox{d}y \ge c\,|x|^{-\ord},\quad
\text{for all $|x|>2$.}
\end{align}
Thus, integral \eqref{infinita} is infinite when $q\leq
({d}/{\ord})$.

To see \eqref{porabajo}, let $x\in \RR^d$ with $|x|>2$. If
$|x-y|<1/{|x|}$, then $|x+y|<\frac{5}{2}|x|$, and
\begin{equation*}
\hbox{e}^{-\frac{1}{4}\left(s|x+y|^2+\frac{1}{s}|x-y|^2\right)}
\ge \hbox{e}^{-2\left(s|x|^2+\frac{1}{s}|x-y|^2\right)}.
\end{equation*}
Therefore, from \eqref{Z_0y1med} and a change of variables, we
have
\begin{align*}
\int_{\RR^d} K_{\ord}(x,y)\,\hbox{d}y &\ge
\int_{\{|x-y|<\frac1{|x|}\}} K_{\ord,0}(x,y)\,\hbox{d}y\\[.3pc]
&\ge \frac{1}{C_1}\int_{\{|x-y|<\frac{1}{|x|}\}}\int_{0}^{{1}/{2}}
s^{-\frac{d}{2}+\ord}\,
\hbox{e}^{-2\left(s|x|^2+\frac{1}{s}|x-y|^2\right)}\, \frac{\hbox{d}s}{s}\, \hbox{d}y\\[.3pc]
&= C\,\int_0^{1/|x|} \int_{0}^{{1}/{2}} s^{-\frac{d}{2}+\ord}\,
\hbox{e}^{-2\left(s|x|^2+\frac{r^2}{s}\right)}\,
\frac{\hbox{d}s}{s}\, r^d\frac{\hbox{d}r}{r}\\[.3pc]
&= C\ |x|^{-2\ord} \int_0^{1} \left(\int_{0}^{{|x|^2}/{2}}
u^{-\frac{d}{2}+\ord}\,
\hbox{e}^{-2\left(u+\frac{t^2}{u}\right)}\,
\frac{\hbox{d}u}{u}\right)\, t^d\,\frac{\hbox{d}t}{t}\\[.3pc]
&\ge C\ |x|^{-2\ord} \int_0^{1} \left(\int_{0}^{1}
u^{-\frac{d}{2}+\ord}\,
\hbox{e}^{-2\left(u+\frac{t^2}{u}\right)}\,
\frac{\hbox{d}u}{u}\right)\, t^d\,\frac{\hbox{d}t}{t}.
\end{align*}
Thus, we have proved \eqref{porabajo}.

In order to show (iv), let $f\in L^p(\RR^d)$ with $1\leq p<
\frac{d}{d-\al}$. By Fubini's theorem, \eqref{laPhia} and
H\"older's inequality,
\begin{align*}
\int_{\RR^d} |H^{-\ord/2} f| &\leq \int_{\{|x|\leq
1\}}\int_{\RR^d}K_{\ord/2}(y,x)\,\hbox{d}y\,|f(x)|\,\hbox{d}x\\[.3pc]
&\quad\, +  \int_{\{|x|> 1\}} \int_{\RR^d}K_{\ord/2}(y,x)\,\hbox{d}y\,|f(x)|\,\hbox{d}x\\[.3pc]
&\leq C\,\|f\|_p +  \|f\|_p \int_{\{|x|>
1\}}\left(\int_{\RR^d}K_{\ord/2}(y,x)\,\hbox{d}y\right)^{p'}
\,\hbox{d}x.
\end{align*}
By the symmetry of $K_{\ord/2}$, inequality \eqref{inteny} and
$p<{d}/({d-\ord})$, the last integral is finite. On the other
hand, if $p\ge {d}/({d-\ord})$ and
\begin{equation*}
f(y)=
\begin{cases}
|y|^{-d+\ord}(\log |y|)^{-1} &\text{for $|y|>2$,}\\[.3pc]
0 & \text{ otherwise,}
\end{cases}
\end{equation*}
then $f$ belongs to $L^p(\RR^d)$ but, from \eqref{porabajo},
\begin{align*}
\int_{\RR^d} |H^{-\ord/2} f| &\ge \int_{\{|y|> 1\}}
\int_{\RR^d}K_{\ord/2}(x,y)\,\hbox{d}x\,|f(y)|\,\hbox{d}y\\[.3pc]
&\ge \int_{\{|y|> 1\}} |y|^{-d}(\log|y|)^{-1} \,\hbox{d}y =
\infty.
\end{align*}

Finally, (v) is a consequence of Remark~\ref{rem:acotH_int}, (i),
(ii), (iii), (iv) and the Riesz--Torin interpolation
theorem.\hfill $\Box$
\end{proof}

As a consequence of the last theorem we have the following
Poincar\'e-type inequalities.

\begin{theor}[\!]\label{teo:poinca} Let $d >1$.
Define the Hermite gradient as
\begin{equation*}
\nabla_H f = (A_{-d}f,\dots,A_{-1}f,A_1f,\dots,A_df).
\end{equation*}
Let $p,q$ in the range $1< p,\ q< \infty$ such that
$\frac{1}{p}-\frac{1}{d} \leq \frac{1}{q}<
\frac{1}{p}+\frac{1}{d}$. Then
\begin{equation*}
\|f\|_q \leq C\, \|\nabla_H f\|_{L^p_{\RR^{2d}}},
\end{equation*}
for all $f\in \LL^p_1$.
\end{theor}

\begin{proof}
It is enough to prove the result for $f\in\D$. From inequality
\eqref{Reves}, we see that
\begin{equation*}
\|f\|_q \leq C\, \sum_{1\leq |j| \leq d} \|R_{j}(f)\|_q,
\end{equation*}
so that, by Lemma~\ref{lem:conmut},
\begin{equation*}
R_jf= \left(\frac{H}{H+2}\right)^{1/2}\,H^{-1/2}A_jf ,
\end{equation*}
where $1\leq j \leq d$.

We have already seen in the proof of Theorem~\ref{teo:acotAj},
that the operator $({H}/({H+2}))^{1/2}$ is bounded on
$L^q(\RR^d)$. Hence, by using Theorem~\ref{teo:acotHal}(v),
\begin{equation*}
\|R_jf\|_q \leq C\,\|A_jf\|_p,
\end{equation*}
where $\frac{1}{p}-\frac{1}{d} \leq \frac{1}{q}<
\frac{1}{p}+\frac{1}{d}$.\hfill $\Box$
\end{proof}

\section{Some applications to Schr\"odinger solutions}\label{sec:schrodinger}

In this section we deal with the unidimensional Schr\"{o}dinger
equation
\begin{equation}\label{scho_eq}
\begin{cases}
\displaystyle i \frac{\partial u(x,t)}{\partial t} = H u(x,t) &x, t \in \RR\\[.5pc]
u(x,0)=f(x)
\end{cases}
\end{equation}
for some initial data $f$.

We are interested in where we have to pick the function $f$ in
order to have almost everywhere convergence of the solution
\begin{equation*}
u(x,t)=\hbox{e}^{itH}f(x)
\end{equation*}
of \eqref{scho_eq} to $f$ as $t$ tends to $0$.

In \cite{1} and \cite{3} the problem with the classical Laplacian
is considered. In \cite{2} the problem for a more general operator
is studied. From that work, it can be derived, for $H$ as a
particular case, that if $f$ belongs to $\LL^p_\ord$ with
$\ord>1$, then we have almost everywhere convergence. Next theorem
gives convergence for orders of differentiability greater than
$1/2$.

\begin{theor}[\!]
If $\ord>1/2$ and $f$ belongs to $\LL_{\ord}^2(\RR)${\rm ,} then
$\hbox{\rm e}^{itH}\!f$ converges to $f$ almost everywhere as $t$
tends to $0$.
\end{theor}

\begin{proof}
If $f$ is a finite linear combination of Hermite functions,
\begin{equation*}
\lim_{t\rightarrow 0} \hbox{e}^{itH}f(x) = f(x)
\end{equation*}
everywhere. Since this kind of functions are dense in
$\LL^2_\ord$, it is enough to prove that the maximal function
\begin{equation*}
T^*f=\sup_{t>0} |\hbox{e}^{itH}f|
\end{equation*}
satisfies the inequality
\begin{equation*}
\int_I T^*f \leq C\, \|f\|_{\LL^2_\ord}\text{,}
\end{equation*}
for all compact interval $I$ of the real line not containing the
origin, and $C$ a constant that may depend on the interval $I$ but
not on $f$.

In order to see this property, we will use the following estimate
of Hermite functions that can be found in \cite{11}
(Theorem~8.91.3, p.~236).

If $I$ is a bounded interval and does not contain the origin,
there exist  constants $C$ and $k_0$ such that
\begin{equation}\label{hfest}
|\hf_k(x)|\leq \frac{C}{k^{1/4}}
\end{equation}
for all $x\in I$ and $k\ge k_0$.

Let $f$ be in $\LL^2_\ord$. As $f$ belongs to $L^2(\RR)$ it can be
written as
\begin{equation*}
f(x)=\sum_{k=0}^\infty a_k\,\hf_k(x).
\end{equation*}

By Tonelli's theorem, estimate \eqref{hfest} and H\"older's
inequality, we get
\begin{align*}
\int_I |T^*f(x)|\,\hbox{d}x &\leq \int_I \sup_{t>0}
\left|\sum_{k=0}^\infty
\,a_k\,\hbox{e}^{it(2k+1)}\,\hf_k(x)\right|\,\hbox{d}x\\[.3pc]
&\leq \sum_{k=0}^\infty\, |a_k|\,\int_I |\hf_k(x)|\,\hbox{d}x\\[.3pc]
&\leq C\, \left(C+\sum_{k=k_0}^\infty
\frac{1}{k^{1/2}(2k+1)^{\ord}}\right)^{1/2}
\left(\sum_{k=0}^\infty a_k^2\, (2k+1)^{\ord} \right)^{1/2}\\[.3pc]
&\leq C \left(C+\sum_{k=k_0}^\infty
\frac{1}{k^{1/2+\ord}}\right)^{1/2}\|f\|_{\LL^2_{\ord}}.
\end{align*}

Since $\ord>1/2$, we have $1/2+\ord>1$ and  the last series is
convergent.\hfill $\Box$
\end{proof}

\begin{theor}[\!]
If $\ord<1/4${\rm ,} then there exists a function $f$ in
$\LL^2_\ord(\RR)$ such that
\begin{equation*}
\lim_{t\rightarrow 0}{\rm e}^{-itH}f(x) = \infty
\end{equation*}
for almost every $x\in\RR$.
\end{theor}

\begin{proof}
For $\ord<1/4$, in \cite{3} the authors find an $f$ belonging to
the classical Sobolev space $L^2_\ord$ and compactly supported so
that
\begin{equation}\label{nolimite}
\liminf_{t\rightarrow 0} |\hbox{e}^{-it\Delta}f(x)| =
\infty\text{.}
\end{equation}

Since $f$ is compactly supported, it follows from
Theorem~\ref{teo:HSvsS}(iii) that $f$ belongs to $\LL^2_\ord$.
Then, it is sufficient to compare the kernels of $\hbox{e}^{-it
\lap}$ and $\hbox{e}^{-itH}$ for small values of $t$. In fact,
\begin{equation*}
\hbox{e}^{-it\lap}f(x) = \int_{\RR}\!
W_{it}(x-y)\,f(y)\,\hbox{d}y,
\end{equation*}
with
\begin{equation*}
W_z(x) = \frac{1}{\sqrt{4\pi z}}
\exp\!\left(\frac{-|x|}{4z}\right),\quad z\in \mathbb{C}
\end{equation*}
and
\begin{equation*}
\hbox{e}^{-itH}f(x) =  \int_{\RR}\! G_{it}(x,y)\,f(y)\,\hbox{d}y,
\end{equation*}
where
\begin{align*}
&G_z(x,y)\\[.4pc]
&\quad\, = \frac{1}{\sqrt{2\pi \sinh(2z)}}
\exp\!\left(-\frac{1}{2} |x-y|^2 \coth(2z) - x\cdot
y\tanh(z)\right),\quad z\in \mathbb{C}.
\end{align*}
Then, for a fixed $x\in\RR$, we have
\begin{align*}
&\lim_{t\rightarrow 0} |W_{it}(x-y) - G_{it}(x,y)|\\[.4pc]
&\quad\,= \frac{1}{\sqrt{2\pi}} \lim_{t\rightarrow 0}
\left|\frac{\hbox{e}^{i\frac{|x-y|^2}{2\tan(2t)} - i x \cdot y
\tan(t)}}{\sqrt{\sin(2t)}} -
\frac{\hbox{e}^{i\frac{|x-y|^2}{4t}}}{\sqrt{2t}}\right| = 0
\end{align*}
uniformly for $y$ in a compact subset of $\RR$. Thus, by
\eqref{nolimite} we also have
\begin{equation*}
\liminf_{t\rightarrow 0} |\hbox{e}^{-itH}f(x)| = \infty.
\end{equation*}

$\left.\right.$\vspace{-1.5pc}

\hfill $\Box$
\end{proof}

\section*{Acknowledgement}

The first author is grateful to the Department of Mathematics at
Universidad Aut\'onoma de Madrid for its hospitality during the
period of this research. The authors were supported by MCYT
2002-04013-C02-02 and HARP network HPRN-CT-2001-00273 of the
European Commission, and CONICET (Argentina).

\end{document}